# Multicolor urn models with reducible replacement matrices

ARUP BOSE[*], AMITES DASGUPTA[**] and KRISHANU MAULIK[†]

*Statistics and Mathematics Unit, Indian Statistical Institute, 202 B.T. Road, Kolkata 700108, India. E-mail: [*]abose@isical.ac.in; [**]amites@isical.ac.in; [†]krishanu@isical.ac.in*

Consider the multicolored urn model where, after every draw, balls of the different colors are added to the urn in a proportion determined by a given stochastic replacement matrix. We consider some special replacement matrices which are not irreducible. For three- and four-color urns, we derive the asymptotic behavior of linear combinations of the number of balls. In particular, we show that certain linear combinations of the balls of different colors have limiting distributions which are variance mixtures of normal distributions. We also obtain almost sure limits in certain cases in contrast to the corresponding irreducible cases, where only weak limits are known.

## 1. Introduction

Consider an urn model with balls of $K$ colors. The row vector $\boldsymbol{C}_0$ will denote the number of balls of each color we start with. (By abuse of terminology, we shall allow the number of balls to be any non-negative real number.) The vector $\boldsymbol{C}_0$ will be taken to be a probability vector, that is, each coordinate is non-negative and the coordinates add up to 1. Suppose $R = ((r_{ij}))$ is a $K \times K$ non-random stochastic (i.e., each row sum is one) replacement matrix. The results of this paper extend to non-random replacement matrices with constant (not necessarily one) row sums by an obvious rescaling. Let $\boldsymbol{C}_n$ be the row vector giving the number of balls of each color after the $n$th trial. At the $n$th trial, a ball is drawn at random, and so a ball of $i$th color appears with probability $\boldsymbol{C}_{i,n-1}/n$. If a ball of $i$th color appears, then the number of balls of $j$th color is increased by $r_{ij}$. If $R$ equals the identity matrix, then it is well known (see, e.g., [3]) that $\boldsymbol{C}_n/(n+1)$ converges almost surely to a Dirichlet random vector with parameters given by the starting vector $\boldsymbol{C}_0$.

Let $\mathbf{1}$ or $\mathbf{0}$ stand respectively for the column vector of relevant dimension with all coordinates 1 or 0. For any vector $\boldsymbol{\xi}$, $\boldsymbol{\xi}^2$ will be the vector whose coordinates are the square of those of $\boldsymbol{\xi}$.

In Section 2, we consider two color models ($K = 2$). If the replacement matrix $R$ is not the identity matrix, then it has two right eigenvectors, $\mathbf{1}$ and $\boldsymbol{\xi}$ corresponding to the principal eigenvalue 1 and the non-principal eigenvalue $\lambda$, respectively, with $|\lambda| < 1$. If $R$ is irreducible, the asymptotic properties of $\boldsymbol{C}_n\mathbf{1}$ and $\boldsymbol{C}_n\boldsymbol{\xi}$ are well known in the literature, see Proposition 2.1.

When the replacement matrix $R$ is reducible but not the identity matrix, then, after possibly interchanging the names of the colors, $R$ is an upper triangular matrix

$$R = \begin{pmatrix} s & 1-s \\ 0 & 1 \end{pmatrix}, \tag{1}$$

for $0 < s < 1$. Here the non-principal eigenvalue is $s$ with the corresponding eigenvector $\boldsymbol{\xi} = (1,0)'$. The asymptotic behavior of the linear combinations is given in Proposition 2.2. In this case, $\boldsymbol{C}_n\boldsymbol{\xi}/n^s = W_n/n^s$ converges almost surely for all values of $s$ in contrast to the irreducible case. See also Theorems 1.3(v), 1.7, 1.8 and 8.8 of [7], where the distribution of the limiting random variable was identified using methods from the branching process.

In the multicolor case, when $R$ is irreducible, the weak/strong laws corresponding to different linear combinations are completely known, see [1, 6]. Gouet [4] considered (reducible) replacement matrices that are block diagonal, with all but the last block irreducible. The last block was taken to be block upper triangular, which cannot be converted into a block diagonal one, and each diagonal subblock of the last block was assumed to be a multiple of some irreducible stochastic matrix. He showed (cf. Theorem 3.1 of [4]) that the proportions of colors converge almost surely to a constant vector where the non-zero coordinates correspond to all but the last diagonal block and the last diagonal subblock of the last diagonal block. We call the corresponding colors *dominant*. To avoid trivial situations, we shall always assume positive contribution to at least one non-dominant color in the initial vector $\boldsymbol{C}_0$.

We shall consider three- and four-color urn models with block upper triangular replacement matrices that are not block diagonal. The diagonal blocks will be taken to be irreducible and we shall extend the result obtained in [4] by obtaining the limiting results for linear combinations corresponding to a complete set of linearly independent vectors.

Specifically, in Section 3, we consider three colors – white, black and green – and the $3 \times 3$ replacement matrix

$$R = \begin{pmatrix} sQ & & 1-s \\ & & 1-s \\ 0 & 0 & 1 \end{pmatrix}, \tag{2}$$

where $0 < s < 1$, and $Q$ is a $2 \times 2$ irreducible aperiodic stochastic matrix with stationary distribution $\boldsymbol{\pi}_Q$. Here green alone is the dominant color and we assume that $W_0 + B_0 > 0$. We show in Theorem 3.1(iv) that $(W_n, B_n)/n^s \overset{\text{a.s.}}{\to} \boldsymbol{\pi}_Q V$, where $P(V > 0) = 1$ and $V$ is non-degenerate. If $\boldsymbol{\xi}$ is the eigenvector corresponding to the non-principal eigenvalue $\lambda$ of $Q$, weak/strong laws for $(W_n, B_n)\boldsymbol{\xi}$ are also provided in Theorem 3.1. If $\lambda \leq 1/2$, then the weak limit is a variance mixture of normal, in contrast to the irreducible model, where the weak limit is normal.

In Section 4, we consider another type of reducible replacement matrix with two dominant colors:

$$R = \begin{pmatrix} s & (1-s)\boldsymbol{p} \\ \boldsymbol{0} & P \end{pmatrix}, \tag{3}$$

where $P$ is a $2 \times 2$ irreducible stochastic matrix, $\boldsymbol{p}$ is a row probability vector and $0 < s < 1$. If the eigenvalues of $P$ are $\lambda$ and 1, then $s$, $\lambda$ and 1 are eigenvalues of $R$. Clearly $(1,0,0)'$ is the eigenvector corresponding to $s$ and the behavior of the corresponding linear combination, $W_n$, follows directly from Proposition 2.2.

Now consider the eigenvalue $\lambda$. If $R$ is diagonalizable, then the weak/strong law of the linear combination given by the eigenvector corresponding to $\lambda$ is summarized in Theorem 4.1. If $R$ is not diagonalizable, then one of the eigenvalues is repeated, namely $\lambda = s$, and the repeated eigenvalue has eigenspace of dimension 1, spanned by $(1,0,0)'$. Consider the Jordan decomposition of $R$, $RT = TJ$, where $T$ is non-singular and

$$J = \begin{pmatrix} s & 1 & 0 \\ 0 & s & 0 \\ 0 & 0 & 1 \end{pmatrix}. \tag{4}$$

The first and the third columns of $T$ can be chosen as $(1,0,0)'$ and $\boldsymbol{1}$ respectively. For the linear combination corresponding to the middle column, we get weak/strong law. The convergence is in the almost sure sense, whenever $\lambda \geq 1/2$, unlike the irreducible and the diagonalizable reducible cases. For $\lambda = 1/2$, in the irreducible and the diagonalizable reducible cases, we have weak convergence only. Also, the scaling for the irreducible case is $\sqrt{n \log^3 n}$ and for the diagonalizable reducible case is $\sqrt{n \log n}$, unlike the non-diagonalizable reducible case, where the scaling is $\sqrt{n \log^2 n}$.

Other interesting reducible three-color urn models have been considered in the literature. For example, [2] and [9] consider three-color urn models with triangular replacement matrices. Our emphasis is on replacement matrices with block triangular structure, given by (2) and (3). Note that $Q$ in (2) is assumed to be a stochastic matrix. However, our techniques do not have a direct extension to the case where $Q$ does not have a constant row sum. In another related work, Pouyanne [8] allows eigenvalues of the replacement matrix to be complex and obtains interesting results for appropriate linear combinations. For example, in his Theorems 3.5 and 3.6, rates are given for the linear combinations corresponding to the eigenvalue with the second largest real part, where it is bigger than $1/2$. In our setup, all the eigenvalues are real and we obtain the rates for all possible linear combinations.

The results for three-color urn models are extended to four-color (white, black, green and yellow) urns with the reducible replacement matrix given by

$$R = \begin{pmatrix} sQ & E \\ 0 & P \end{pmatrix}, \tag{5}$$

where each component is a $2 \times 2$ matrix and furthermore $P$ and $Q$ are irreducible stochastic matrices, $0 < s < 1$. The results are summarized in Propositions 4.2–4.5. An interesting phenomenon is observed in Proposition 4.5, where the replacement matrix is not diagonalizable and the repeated eigenvalue is zero. Unlike the behavior of the corresponding linear combination in other cases, where it remains a constant, we get a weak limit of variance mixture of normal distribution in this case.

Before proceeding with the details, note that the proofs are based on studying the behavior of appropriate martingales with the filtration $\mathcal{F}_n$ being the natural filtration of the sequence $\{\boldsymbol{C}_n\}$.

# 2. Two-color urn models

Define

$$\Pi_n(\lambda) = \prod_{j=0}^{n-1}\left(1 + \frac{\lambda}{j+1}\right). \tag{6}$$

Recall that Euler's formula for gamma function gives

$$\Pi_n(\lambda) \sim n^\lambda/\Gamma(\lambda+1), \qquad \lambda \text{ not a negative integer.} \tag{7}$$

This will be used at several places later.

We first mention the asymptotic behavior in two-color irreducible urn models. The following results are well known. See, for example, [1, 6].

**Proposition 2.1.** *In a two-color urn model with irreducible replacement matrix $R$,*

$$\frac{\boldsymbol{C}_n}{n+1} \overset{\text{a.s.}}{\to} \boldsymbol{\pi}_R, \tag{8}$$

*where $\boldsymbol{\pi}_R$ is the stationary distribution of $R$. Further, we have:*

  (i) *If $\lambda < 1/2$, then $\boldsymbol{C}_n\boldsymbol{\xi}/\sqrt{n} \Rightarrow N(0, \frac{\lambda^2}{1-2\lambda}\boldsymbol{\pi}_R\boldsymbol{\xi}^2)$.*
  (ii) *If $\lambda = 1/2$, then $\boldsymbol{C}_n\boldsymbol{\xi}/\sqrt{n\log n} \Rightarrow N(0, \lambda^2\boldsymbol{\pi}_R\boldsymbol{\xi}^2)$.*
  (iii) *If $\lambda > 1/2$, then $\boldsymbol{C}_n\boldsymbol{\xi}/\Pi_n(\lambda)$ is an $L^2$-bounded martingale and converges almost surely, as well as in $L^2$, to a non-degenerate random variable.*

**Remark 2.1.** Since $W_n + B_n = n+1$, we have from (8),

$$(W_n, B_n)/(W_n + B_n) \overset{\text{a.s.}}{\to} \boldsymbol{\pi}_R. \tag{9}$$

**Remark 2.2.** From (8), we see that $\boldsymbol{C}_n\boldsymbol{\xi}/(n+1) \overset{\text{a.s.}}{\to} \boldsymbol{\pi}_R\boldsymbol{\xi}$. However $\boldsymbol{\pi}_R\boldsymbol{\xi} = \boldsymbol{\pi}_R R\boldsymbol{\xi} = \lambda\boldsymbol{\pi}_R\boldsymbol{\xi}$, and since $\lambda \neq 1$, we have $\boldsymbol{\pi}_R\boldsymbol{\xi} = 0$. This explains the appropriate scaling up of $\boldsymbol{C}_n\boldsymbol{\xi}/(n+1)$ to obtain the weak laws for the proportions above.

**Remark 2.3.** When $\lambda > 1/2$, using (7), $\boldsymbol{C}_n\boldsymbol{\xi}/n^\lambda$ converges almost surely, as well as in $L^2$, to a non-degenerate random variable. Thus, the scalings in Proposition 2.1(i) and

(iii) are different. Also, in (iii), the distribution of the limit random variable depends on the starting value $(W_0, B_0)$ unlike in (i) and (ii). Furthermore, as in Remark 2.1, we can conclude that, when $\lambda > 1/2$, $(W_n, B_n)\boldsymbol{\xi}/(W_n + B_n)^{\lambda}$ converges almost surely, as well as in $L^2$, to a non-degenerate random variable.

*Remark 2.4.* If $\lambda = 0$, both rows of $R$ equal $\boldsymbol{\pi}_R$. Since $\boldsymbol{\pi}_R\boldsymbol{\xi} = 0$ and clearly $\boldsymbol{C}_n = \boldsymbol{C}_0 + n\boldsymbol{\pi}_R$, we have $\boldsymbol{C}_n\boldsymbol{\xi} = \boldsymbol{C}_0\boldsymbol{\xi}$ for all $n$.

Next, we consider the almost sure limit behavior of the two-color urn model with an upper triangular reducible replacement matrix given by (1).

**Proposition 2.2.** *In a two-color urn model with an upper triangular replacement matrix given by (1), we have:*

  (i) $\boldsymbol{C}_n\mathbf{1}/(n+1) = 1$.
  (ii) $\boldsymbol{C}_n/(n+1) \overset{\text{a.s.}}{\to} (0,1)$.
  (iii) $\boldsymbol{C}_n\boldsymbol{\xi}/\Pi_n(s) = W_n/\Pi_n(s)$ *is an $L^2$-bounded martingale, where $\Pi_n(s)$ is given by (6). Further, $W_n/n^s$ converges to a non-degenerate, positive random variable almost surely, as well as in $L^2$.*

**Proof.** Statement (i) is trivial. Statement (ii) is same as that of (8) in Proposition 2.1 and a proof can be obtained from Proposition 4.3 of [4].

For (iii), observe that the number of white balls evolves as

$$W_{n+1} = W_n + s\chi_{n+1},$$

where $\chi_n$ is the indicator of a white ball in $n$th trial. Define the martingale sequence $V_n = W_n/\Pi_n(s)$, $n \geq 1$. We shall show that $\{V_n\}$ is an $L^2$-bounded martingale and hence converges almost surely, as well as in $L^2$. Also, the variance of $V_n$ increases to that of the limit and hence the limit is non-degenerate. The proposition then follows from (7), the distribution of $V$ and the fact that it is almost surely positive, all of which have been established using branching process techniques in Theorem 1.3(v) of [7].

Clearly, we have

$$V_{n+1} - V_n = \frac{s}{\Pi_{n+1}(s)}\left(\chi_{n+1} - \frac{W_n}{n+1}\right),$$

and further, using $V_{n+1} = V_n + (V_{n+1} - V_n)$ and the martingale property, there exists $N$ (non-random), such that for all $n \geq N$,

$$E[V_{n+1}^2|\mathcal{F}_n] = V_n^2 + \frac{s^2}{\Pi_{n+1}^2(s)}\left[\frac{W_n}{n+1} - \frac{W_n^2}{(n+1)^2}\right]$$

$$\leq V_n^2 + \frac{V_n}{(n+1)\Pi_n(s)}$$

$$\leq V_n^2 + \Gamma(s+1)\frac{1 + V_n^2}{(n+1)^{s+1}} = V_n^2\left[1 + \frac{\Gamma(s+1)}{(n+1)^{s+1}}\right] + \frac{\Gamma(s+1)}{(n+1)^{s+1}}.$$

The last inequality holds for $n \geq N$ and follows from the fact that $V_n \leq (1 + V_n^2)/2$ and (7). Taking further expectation and adding 1 to both sides, we have, for $n \geq N$,

$$E[V_{n+1}^2] + 1 \leq \left[1 + \frac{\Gamma(s+1)}{(n+1)^{s+1}}\right](E[V_n^2] + 1).$$

Iterating, we get for $n \geq N$,

$$E[V_{n+1}^2] + 1 \leq (E[V_N^2] + 1) \prod_{j=N}^{n} \left[1 + \frac{\Gamma(s+1)}{(j+1)^{s+1}}\right]$$

and since $s > 0$, we further have for all $n > N$,

$$E[V_n^2] \leq (E[V_N^2] + 1) \exp\left(\Gamma(s+1) \sum_0^{\infty} j^{-(s+1)}\right) < \infty,$$

which shows $\{V_n\}$ is $L^2$-bounded as required. $\qquad \square$

## 3. One dominant color, $K = 3$

Now we are ready to consider the three-color urn model with only one dominant color, say green. We shall denote the row subvector corresponding to the non-dominant colors $(W_n, B_n)$ as $\boldsymbol{S}_n$. We collect the results in the following theorem.

**Theorem 3.1.** *Consider a three-color urn model with a reducible replacement matrix $R$ given by (2). Suppose the non-principal eigenvalue of $Q$ is $\lambda$ and the corresponding eigenvector is $\boldsymbol{\xi}$. Then the following hold:*

(i) *$\boldsymbol{C}_n \mathbf{1}/(n+1) = 1$.*

(ii) *$\boldsymbol{C}_n/(n+1) \overset{\text{a.s.}}{\to} (0, 0, 1)$.*

(iii) *$\boldsymbol{S}_n \mathbf{1}/(n+1)^s$ converges almost surely, as well as in $L^2$, to a non-degenerate positive random variable $U$.*

(iv) *$\boldsymbol{S}_n/(n+1)^s \overset{\text{a.s.}}{\to} \boldsymbol{\pi}_Q U$.*

(v) *If $\lambda < 1/2$, then $\boldsymbol{S}_n \boldsymbol{\xi}/n^{s/2} \Rightarrow N(0, \frac{s^2 \lambda^2}{s(1-2\lambda)} U \boldsymbol{\pi}_Q \boldsymbol{\xi}^2)$.*

(vi) *If $\lambda = 1/2$, then $\boldsymbol{S}_n \boldsymbol{\xi}/\sqrt{n^s \log n} \Rightarrow N(0, \ s^2 \lambda^2 U \boldsymbol{\pi}_Q \boldsymbol{\xi}^2)$.*

(vii) *If $\lambda > 1/2$, then $\boldsymbol{S}_n \boldsymbol{\xi}/\Pi_n(s\lambda)$ is an $L^2$-bounded martingale and almost surely, as well as in $L^2$, $\boldsymbol{S}_n \boldsymbol{\xi}/n^{s\lambda} \to V$, where $V$ is a non-degenerate random variable.*

*The random variable $U$ in* (iv), (v) *and* (vi) *is the same limiting random variable obtained in* (iii). *The distributions of $U$ and $V$ depend on the initial value $\boldsymbol{S}_0$.*

**Remark 3.1.** Note that the eigenvalues of $R$ are 1, $s$ and $s\lambda$ with corresponding eigenvectors $\mathbf{1}$, $(1, 1, 0)'$ and $(\boldsymbol{\xi}', 0)'$, respectively, yielding the linear combinations $\boldsymbol{C}_n \mathbf{1}$, $\boldsymbol{S}_n \mathbf{1}$ and $\boldsymbol{S}_n \boldsymbol{\xi}$.

**Proof of Theorem 3.1.** Statement (i) is immediate. Statement (ii) follows from Theorem 3.1 and Proposition 4.3 of [4].

Note that $\boldsymbol{S}_n \mathbf{1} = W_n + B_n$. From the structure of $R$, the pair $(W_n + B_n, G_n)$ yields a two-color model with *reducible* replacement matrix

$$\begin{pmatrix} s & 1-s \\ 0 & 1 \end{pmatrix}.$$

Statement (iii) then follows from Proposition 2.2. The distribution of $U$ has been identified in Theorem 1.3(v) of [7].

Consider the successive times $\tau_k$ when either a white or a black ball is observed. Due to the assumed special structure of the matrix $R$, it is only at these times that more white or black balls are added and the total number added is the constant $s$. Thus $\boldsymbol{S}_{\tau_k}/\boldsymbol{S}_{\tau_k}\mathbf{1}$ are the proportions from the evolution of a two-color urn model governed by the irreducible replacement matrix $Q$. Hence by the two-color urn result (9), it converges almost surely to $\boldsymbol{\pi}_Q$. Note that at all other $n$, $\tau_k < n < \tau_{k+1}$, the vector $\boldsymbol{S}_n = \boldsymbol{S}_{\tau_k}$ and hence the ratio is unchanged. Moreover, from the statement (iii), we have $\boldsymbol{S}_n\mathbf{1}/n^s \stackrel{\text{a.s.}}{\to} V$. Now combining all of the above, the proof of the statement (iv) is complete.

For (vii), let $\boldsymbol{\chi}_n$ be the row vector, which takes values $\boldsymbol{\chi}_n = (1,0)$, $(0,1)$ or $(0,0)$ accordingly as the white, black or green balls are observed in $n$th trial and consider the martingale $T_n = S_n\boldsymbol{\xi}/\Pi_n(s\lambda)$. Then the martingale difference is

$$T_{n+1} - T_n = \frac{s\lambda}{\Pi_{n+1}(s\lambda)}\left(\boldsymbol{\chi}_n\boldsymbol{\xi} - \frac{\boldsymbol{S}_n\boldsymbol{\xi}}{n+1}\right)$$

and hence

$$E[T_{n+1}^2] = E[T_n^2] + \left(\frac{s\lambda}{\Pi_{n+1}(s\lambda)}\right)^2 E\left[\frac{S_n\boldsymbol{\xi}^2}{n+1} - \left(\frac{S_n\boldsymbol{\xi}}{n+1}\right)^2\right]$$

$$= E[T_n^2]\left[1 - \frac{(s\lambda)^2}{(n+1)^2(1+s\lambda/(n+1))^2}\right] + \left(\frac{s\lambda}{\Pi_{n+1}(s\lambda)}\right)^2 \frac{1}{(n+1)^{1-s}} E\left[\frac{S_n\boldsymbol{\xi}^2}{(n+1)^s}\right].$$

The first term is bounded by $E[T_n^2]$. From statement (iii), $S_n\mathbf{1}/(n+1)^s$ is $L^2$-bounded and hence $L^1$-bounded. So $S_n\boldsymbol{\xi}^2/(n+1)^s$ is also $L^1$-bounded. Thus, using (7), the second term is bounded by a constant multiple of $n^{-(1+s(2\lambda-1))}$, which is summable as $\lambda > 1/2$. Thus $\{T_n\}$ is an $L^2$-bounded martingale and hence converges almost surely as well as in $L^2$.

For (v) and (vi), we start with the case $\lambda < 1/2$. Call $X_n = \boldsymbol{S}_n\boldsymbol{\xi}/n^{s/2}$. We have the evolution equation for $\boldsymbol{S}_n\boldsymbol{\xi}$ given by

$$\boldsymbol{S}_{n+1}\boldsymbol{\xi} = \boldsymbol{S}_n\boldsymbol{\xi} + s\boldsymbol{\chi}_{n+1}Q\boldsymbol{\xi} = \boldsymbol{S}_n\boldsymbol{\xi} + \lambda s\boldsymbol{\chi}_{n+1}\boldsymbol{\xi}. \tag{10}$$

We now use the decomposition of $X_{n+1}$ into a conditional expectation and a martingale difference

$$X_{n+1} = E(X_{n+1}|\mathcal{F}_n) + \{X_{n+1} - E(X_{n+1}|\mathcal{F}_n)\}.$$

Using (10) and the fact that $(1 + 1/n)^{-s/2} = (1 - s/2n) + O(1/n^2)$ we then get

$$
\begin{aligned}
E(X_{n+1}|\mathcal{F}_n) &= \frac{\boldsymbol{S}_n\boldsymbol{\xi}}{n^{s/2}}(1 + 1/n)^{-s/2} + \frac{\lambda s}{(n+1)^{s/2}}\frac{\boldsymbol{S}_n\boldsymbol{\xi}}{n+1} \\
&= X_n\left(1 - \frac{s}{2n} + O\left(\frac{1}{n^2}\right)\right) + \lambda s X_n\left(1 + \frac{1}{n}\right)^{-s/2}\frac{1}{n+1} \\
&= X_n\left(1 - \frac{s(1/2 - \lambda)}{n}\right) + X_n O(n^{-2}).
\end{aligned}
\tag{11}
$$

On the other hand, the martingale difference is really

$$
M_{n+1} := X_{n+1} - E(X_{n+1}|\mathcal{F}_n) = \frac{\lambda s}{(n+1)^{s/2}}\left(\boldsymbol{\chi}_{n+1} - \frac{\boldsymbol{S}_n}{n+1}\right)\boldsymbol{\xi},
\tag{12}
$$

so that

$$
X_{n+1} = X_n\left(1 - \frac{s(1/2 - \lambda)}{n}\right) + X_n O(n^{-2}) + M_{n+1}.
\tag{13}
$$

Iterating the equation above, we get

$$
\begin{aligned}
X_{n+1} &= X_1\prod_{i=1}^{n}\left(1 - \frac{s(1/2 - \lambda)}{i}\right) + \sum_{j=1}^{n}X_j O(j^{-2})\prod_{i=j+1}^{n}\left(1 - \frac{s(1/2 - \lambda)}{i}\right) \\
&\quad + \sum_{j=1}^{n}M_{j+1}\prod_{i=j+1}^{n}\left(1 - \frac{s(1/2 - \lambda)}{i}\right).
\end{aligned}
\tag{14}
$$

Since $\lambda < 1/2$, we have $\Pi_n(-s(1/2 - \lambda)) \sim n^{-s(1/2-\lambda)}/\Gamma(1 - s(1/2 - \lambda)) \to 0$, and hence the first term above converges to 0 for every sample point. The continued product in the second term is bounded by 1. Since the coordinates of $\boldsymbol{S}_n/(n+1)$ are bounded by 1, we have that $|X_n|/n^{1-s/2}$ is bounded for every sample point. Thus the sum of the elements of the second term above is bounded by a multiple of $\sum_1^{\infty}j^{-(1+s/2)}$, which is finite; and individually each element tends to zero since each infinite product diverges to zero. Hence the second term of (14) tends to zero for every sample point.

Now we turn to the third term of (14),

$$
Z_{n+1} = \sum_{j=1}^{n}M_{j+1}\prod_{i=j+1}^{n}\left(1 - \frac{s(1/2 - \lambda)}{i}\right).
\tag{15}
$$

We verify the conditional Lyapunov condition and compute the conditional variance as $n \to \infty$. The conditional Lyapunov condition demands that for some $k > 2$,

$$
\sum_{j=1}^{n}E(|M_{j+1}|^k|\mathcal{F}_j)\prod_{i=j+1}^{n}\left(1 - \frac{s(1/2 - \lambda)}{i}\right)^k \overset{\text{a.s.}}{\to} 0.
$$

Since each coordinate of $\boldsymbol{\chi}_{n+1}$ and $\boldsymbol{S}_n/(n+1)$ is bounded by 1, the martingale difference defined in (12) is bounded by a constant multiple of $(n+1)^{-s/2}$. Thus the above sum is bounded by a constant multiple of

$$\sum_{j=1}^{n} j^{-ks/2} \prod_{i=j+1}^{n} \left(1 - \frac{s(1/2-\lambda)}{i}\right)^k,$$

which tends to zero by the bounded convergence theorem, provided we choose $k > 2/s$.

Now we compute the conditional variance. An exact computation and the statement (iv) yields, with probability 1,

$$E(M_{n+1}^2|\mathcal{F}_n) = \frac{(\lambda s)^2}{(n+1)^s} \left[\frac{\boldsymbol{S}_n \boldsymbol{\xi}^2}{n+1} - \left(\frac{\boldsymbol{S}_n \boldsymbol{\xi}}{n+1}\right)^2\right] \sim \frac{(\lambda s)^2}{n+1} U \boldsymbol{\pi}_Q \boldsymbol{\xi}^2.$$

Then, writing $\prod_{i=j+1}^{n}(1 - \frac{s(1/2-\lambda)}{i}) = \Pi_n(-s(1/2-\lambda))/\Pi_j(-s(1/2-\lambda))$ and using (7), the sum of the conditional variances satisfies, on a set of probability 1,

$$\sum_{j=1}^{n} E(M_{j+1}^2|\mathcal{F}_j) \prod_{i=j+1}^{n} \left(1 - \frac{s(1/2-\lambda)}{i}\right)^2 \sim \frac{(\lambda s)^2 U \boldsymbol{\pi}_Q \boldsymbol{\xi}^2}{n^{s(1-2\lambda)}} \sum_{j=1}^{n} \frac{1}{j^{1-s(1-2\lambda)}},$$

which converges almost surely to $(\lambda s)^2 U \boldsymbol{\pi}_Q \boldsymbol{\xi}^2/s(1-2\lambda)$. Thus, by martingale central limit theorem (see Corollary 3.1 of [5]), the limiting distribution of $Z_{n+1}$, and hence $X_{n+1}$ is the required variance mixture of normal.

Since the analysis for the statement (vi) is similar, we omit the details and provide only a brief sketch of the arguments. We start with $X_n = \boldsymbol{S}_n \boldsymbol{\xi}/\sqrt{n^s \log n}$. The following is the relevant martingale decomposition now. To express the decomposition, the following straightforward approximations are used:

$$(1+1/n)^{-s/2} = 1 - \frac{s}{2n} + \mathrm{O}(n^{-2}) \quad \text{and} \quad \frac{\log n}{\log(n+1)} = \frac{\log n}{\log n + 1/n + \mathrm{O}(n^{-2})}$$

together give

$$\left(\frac{n}{n+1}\right)^{s/2} \sqrt{\frac{\log n}{\log(n+1)}} = \left(1 - \frac{s}{2n} + \mathrm{O}(n^{-2})\right)\left(1 - \frac{1}{2n \log n} + \mathrm{O}\left(\frac{1}{n^2 \log n}\right)\right)$$

$$= 1 - \frac{s}{2n} - \frac{1}{2n \log n} + \mathrm{O}(n^{-2}).$$

Using $\lambda = 1/2$ carefully, the conditional expectation becomes

$$E(X_{n+1}|\mathcal{F}_n) = X_n[1 - (2n \log n)^{-1}] + X_n \mathrm{O}(n^{-2})$$

and, for the martingale difference, we get,

$$M_{n+1} := X_{n+1} - E(X_{n+1}|\mathcal{F}_n) = \frac{s}{2\sqrt{(n+1)^s \log(n+1)}}\left(\boldsymbol{\chi}_{n+1} - \frac{\boldsymbol{S}_n}{n+1}\right)\boldsymbol{\xi}.$$

These together give us the recursion on $X_n$ as

$$X_{n+1} = X_n[1 - (2n\log n)^{-1}] + X_n\mathrm{O}(n^{-2}) + M_{n+1},$$

a decomposition similar to (13). The rest of the proof follows as before with appropriate changes.                                                                                                   □

**Remark 3.2.** Theorem 3.1 gives the scaling for all the linear combinations except when $\lambda = 0$, in which case (v) applies and we obtain $\boldsymbol{S}_n\boldsymbol{\xi}/n^{s/2} \xrightarrow{\mathrm{P}} 0$. However, as discussed in Remark 2.4, $Q$ has both rows the same as $\boldsymbol{\pi}_Q$, which satisfies $\boldsymbol{\pi}_Q\boldsymbol{\xi} = 0$. Since $\boldsymbol{S}_n$ changes only when a white or black ball appears, we have $\boldsymbol{S}_n = \boldsymbol{S}_0 + (W_n + B_n)\boldsymbol{\pi}_Q$ and hence $\boldsymbol{S}_n\boldsymbol{\xi} = \boldsymbol{S}_0\boldsymbol{\xi}$ for all $n$.

# 4. Two dominant colors, $K = 3, 4$

We now consider the three-color case with two dominant colors. The replacement matrix $R$, given by (3), is

$$R = \begin{pmatrix} s & (1-s)\boldsymbol{p} \\ 0 & P \end{pmatrix},$$

where $\boldsymbol{p}$ is a probability vector and $P$ is a $2 \times 2$ irreducible stochastic matrix. Thus 1 is always an eigenvalue of $P$ with the corresponding eigenvector $\mathbf{1}$. We shall denote the other eigenvalue of $P$ as $\lambda$ with corresponding eigenvector $\boldsymbol{\xi}$. Then $s$ and 1 are two eigenvalues of $R$ with corresponding eigenvectors $(1, 0, 0)'$ and $\mathbf{1}$, respectively. Observe that $\boldsymbol{C}_n(1, 0, 0)' = W_n$. The results for this linear combination follow from two-color urn model results, and we summarize them below. We shall denote the stationary distribution of $P$ by $\boldsymbol{\pi}_P$.

**Proposition 4.1.** *Consider a three-color urn model with two dominant colors and the replacement matrix given by (3). Then:*

   (i) $\boldsymbol{C}_n\mathbf{1}/(n+1) = 1$.
   (ii) $\boldsymbol{C}_n/n \xrightarrow{\mathrm{a.s.}} (0, \boldsymbol{\pi}_P)$.
   (iii) $W_n/n^s \to V$ *almost surely, as well as in* $L^2$.

*In* (iii), *if we start with the initial vector* $\boldsymbol{C}_0 = (W_0, B_0, G_0)$, *then* $V$ *has the same distribution as the limit random variable in Theorem* 3.1(vii) *with the initial vector* $(W_0, B_0 + G_0)$.

**Proof.** Statement (i) is trivial. The proof of (ii) is given in Theorem 3.1 or Proposition 4.3 of [4]. For the remaining part, consider the two-color urn model $(W_n, B_n + G_n)$ obtained by collapsing the last two colors. This will have the replacement matrix as in (1) and the results will follow from Proposition 2.2.                                                                       □

However, the one remaining linear combination is more subtle. The choice of the linear combination depends on whether $R$ is similar to a diagonal matrix or, equivalently, has a complete set of eigenvectors. Suppose $\lambda \neq s$. Then $R$ is diagonalizable and $\boldsymbol{v}_2 = (c, \boldsymbol{\xi}')'$ is an eigenvector of $R$ corresponding to $\lambda$, with $c = (1-s)\boldsymbol{p}\boldsymbol{\xi}/(\lambda - s)$. If $\lambda = s$, then $R$ is diagonalizable if and only if $\boldsymbol{p}\boldsymbol{\xi} = 0$. In that case, $(0, \boldsymbol{\xi}')'$ is another eigenvector of $R$ corresponding to $s$ independent of $(1, 0, 0)'$. Also note that, in that case, $\boldsymbol{p}$ is orthogonal to $\boldsymbol{\xi}$ and since $\boldsymbol{p}$ is a probability vector we have $\boldsymbol{p} = \boldsymbol{\pi}_P$. In the diagonalizable case with $\lambda = s$, we denote this remaining vector $(0, \boldsymbol{\xi}')'$ by $\boldsymbol{v}_2$ and consider the corresponding linear combination. The following theorem summarizes the results.

**Theorem 4.1.** *Consider a three-color urn model with replacement matrix $R$ given by (5), where $R$ is diagonalizable. Then the following weak/strong laws hold:*

(i) *If $\lambda < 1/2$, then $\boldsymbol{C}_n \boldsymbol{v}_2/\sqrt{n} \Rightarrow N(0, \frac{\lambda^2}{1-2\lambda}\boldsymbol{\pi}_P \boldsymbol{\xi}^2)$.*

(ii) *If $\lambda = 1/2$, then $\boldsymbol{C}_n \boldsymbol{v}_2/\sqrt{n \log n} \Rightarrow N(0, \lambda^2 \boldsymbol{\pi}_P \boldsymbol{\xi}^2)$.*

(iii) *If $\lambda > 1/2$, then $\boldsymbol{C}_n \boldsymbol{v}_2/\Pi_n(\lambda)$ is an $L^2$-bounded martingale and $\boldsymbol{C}_n \boldsymbol{v}_2/n^\lambda$ converges almost surely to a non-degenerate random variable.*

**Proof.** The proofs of (i) and (ii) are similar to those of (v) and (vi) of Theorem 3.1, so we omit them.

Define $\boldsymbol{\chi}_n$ as the row vector that takes values $(1, 0, 0)$, $(0, 1, 0)$ and $(0, 0, 1)$ accordingly as white, black or green balls appear in $n$th trial. Also define $Z_n = \boldsymbol{C}_n \boldsymbol{v}_2/\Pi_n(\lambda)$. It is simple to check that $\{Z_n\}$ is a martingale. Note that,

$$Z_{n+1} - Z_n = \frac{\lambda}{\Pi_{n+1}(\lambda)}\left(\boldsymbol{\chi}_{n+1} - \frac{\boldsymbol{C}_n}{n+1}\right)\boldsymbol{v}_2,$$

which gives us

$$E[(Z_{n+1} - Z_n)^2 | \mathcal{F}_n] = \frac{\lambda^2}{\Pi_{n+1}^2(\lambda)}\left[\frac{\boldsymbol{C}_n \boldsymbol{v}_2^2}{n+1} - \left(\frac{\boldsymbol{C}_n \boldsymbol{v}_2}{n+1}\right)^2\right].$$

Also, $\boldsymbol{C}_n/(n+1)$ is bounded by 1 for each coordinate. Hence, the conditional expectation above is bounded by a constant multiple of $n^{-2\lambda}$. So we get $E[Z_{n+1}^2] = \sum_{i=1}^n \{E[(Z_{i+1} - Z_i)^2]\}$ is bounded by a constant multiple of $\sum_1^\infty i^{-2\lambda}$, which is finite, as $\lambda > 1/2$. Thus, $\{Z_n\}$ is $L^2$-bounded. The rest of the statement (iii) follows from (7). $\square$

If $R$ is not diagonalizable, then a complete set of eigenvectors is not available and one of the eigenvalues must be repeated, which gives $s = \lambda$ and $\boldsymbol{p} \neq \boldsymbol{\pi}_P$. So we consider the Jordan decomposition $RT = TJ$, where $J$ is given by (4). We can choose the first and third columns of $T$ as $\boldsymbol{t}_1 = (1, 0, 0)'$ and $\boldsymbol{t}_3 = \boldsymbol{1}$. Also the subvector of the lower two coordinates of $\boldsymbol{t}_2$ is an eigenvector of $P$ corresponding to $s$. We shall denote it by $\boldsymbol{\xi}$ as well. The behavior of $\boldsymbol{C}_n \boldsymbol{t}_2$ is substantially different from the irreducible case given in Theorem 3.15 of [6] or the diagonalizable case in Theorem 4.1 above.

**Theorem 4.2.** *Consider a three-color urn model with replacement matrix $R$ given by (5), where $R$ is not diagonalizable. Then, we have:*

(i) *If $s < 1/2$, then $\boldsymbol{C}_n \boldsymbol{t}_2/\sqrt{n} \Rightarrow N(0, \frac{s^2}{1-2s} \pi_P \boldsymbol{\xi}^2)$.*

(ii) *If $s \geq 1/2$, then $\boldsymbol{C}_n \boldsymbol{t}_2/n^s \log n$ converges to $V$ almost surely, as well as in $L^2$, where $V$ is the almost sure limit random variable obtained in Proposition 4.1(iii).*

**Proof.** We first consider the case when $s < 1/2$. Call $X_n = \boldsymbol{C}_n \boldsymbol{t}_2/\sqrt{n}$. Define the row vector $\boldsymbol{\chi}_n$ as in the proof of Theorem 4.1. We shall split $X_{n+1}$ into conditional expectation and martingale difference parts as in the proof of Theorem 3.1(v). From the Jordan decomposition of $R$ and the form (4) of $J$, the evolution equation for $\boldsymbol{C}_n$ is given by

$$\boldsymbol{C}_{n+1} \boldsymbol{t}_2 = \boldsymbol{C}_n \boldsymbol{t}_2 + s \boldsymbol{\chi}_{n+1} \boldsymbol{t}_2 + \boldsymbol{\chi}_{n+1} \boldsymbol{t}_1.$$

Hence the conditional expectation becomes

$$E(X_{n+1} | \mathcal{F}_n) = \frac{\boldsymbol{C}_n \boldsymbol{t}_2}{\sqrt{n+1}} \left( 1 + \frac{s}{n+1} \right) + \frac{1}{(n+1)^{3/2}} \boldsymbol{C}_n \boldsymbol{t}_1$$

$$= X_n \left( 1 - \frac{1/2 - s}{n+1} \right) + X_n O(n^{-2}) + \frac{1}{(n+1)^{3/2}} W_n,$$

since $\boldsymbol{C}_n \boldsymbol{t}_1 = W_n$. Using the notation $s\boldsymbol{t} = \boldsymbol{t}_1 + s\boldsymbol{t}_2$, the martingale difference term becomes

$$M_{n+1} := X_{n+1} - E(X_{n+1} | \mathcal{F}_n) = \frac{s}{\sqrt{n+1}} \left( \boldsymbol{\chi}_{n+1} - \frac{\boldsymbol{C}_n}{n+1} \right) \boldsymbol{t}.$$

Putting this together, we get a recursion on $X_n$ as

$$X_{n+1} = X_n \left( 1 - \frac{1/2 - s}{n} \right) + X_n O(n^{-2}) + \frac{W_n}{(n+1)^{3/2}} + M_{n+1},$$

and iterating we get,

$$X_{n+1} = X_1 \prod_{i=1}^{n} \left( 1 - \frac{1/2 - s}{i} \right) + \sum_{j=1}^{n} X_j O(j^{-2}) \prod_{i=j+1}^{n} \left( 1 - \frac{1/2 - s}{i} \right)$$

$$+ \sum_{j=1}^{n} \frac{W_j}{(j+1)^{3/2}} \prod_{i=j+1}^{n} \left( 1 - \frac{1/2 - s}{i} \right) + \sum_{j=1}^{n} M_{j+1} \prod_{i=j+1}^{n} \left( 1 - \frac{1/2 - s}{i} \right),$$

which is similar to the decomposition (14), except for the additional third term. Further analysis is similar to that done for Theorem 3.1(v), except for the contribution of the third term, which we now show to be negligible with probability 1. By Proposition 4.1(iii), writing $\prod_{i=j+1}^{n} (1 - \frac{1/2-s}{i}) = \Pi_n(-(1/2-s))/\Pi_j(-(1/2-s))$ and using (7), the third term

is of the order of

$$\frac{1}{n^{1/2-s}}\sum_{j=1}^{n}\frac{V}{j^{3/2-s}}\frac{1}{j^{-(1/2-s)}} \sim \frac{V\log n}{n^{1/2-s}} \to 0$$

almost surely, since $s < 1/2$.

Using Proposition 4.1(ii), the structure of the vectors $\boldsymbol{t}_1$ and $\boldsymbol{t}_2$ and the fact $\boldsymbol{\pi}_P\boldsymbol{\xi} = 0$, the conditional variance term is

$$E(M_{n+1}^2|\mathcal{F}_n) = \frac{s^2}{n+1}\left[\frac{\boldsymbol{C}_n\boldsymbol{t}^2}{n+1} - \left(\frac{\boldsymbol{C}_n\boldsymbol{t}}{n+1}\right)^2\right] \sim \frac{s^2}{n+1}\boldsymbol{\pi}_P\boldsymbol{\xi}^2,$$

which gives the required variance for the limiting normal distribution.

Now we consider the other situation, where $s \geq 1/2$. Using the form (4) of $J$ in the Jordan decomposition of $R$, we again have $R\boldsymbol{t}_2 = \boldsymbol{t}_1 + s\boldsymbol{t}_2 = s\boldsymbol{t}$. Thus

$$\boldsymbol{C}_{n+1}\boldsymbol{t}_2 = \boldsymbol{C}_n\boldsymbol{t}_2 + \boldsymbol{\chi}_{n+1}R\boldsymbol{t}_2 = \boldsymbol{C}_n\boldsymbol{t}_2 + s\boldsymbol{\chi}_{n+1}\boldsymbol{t},$$

which implies

$$E[\boldsymbol{C}_{n+1}\boldsymbol{t}_2|\mathcal{F}_n] = \boldsymbol{C}_n\boldsymbol{t}_2\left(1 + \frac{s}{n+1}\right) + \frac{\boldsymbol{C}_n}{n+1}\boldsymbol{t}_1.$$

This gives us the martingale

$$X_n = \frac{\boldsymbol{C}_n\boldsymbol{t}_2}{\Pi_n(s)} - \sum_{j=1}^{n-1}\frac{1}{j+1}\frac{\boldsymbol{C}_j\boldsymbol{t}_1}{\Pi_{j+1}(s)}. \tag{16}$$

The martingale difference is then given by $X_{n+1} - X_n = s(\boldsymbol{\chi}_{n+1} - \frac{\boldsymbol{C}_n}{n+1})\boldsymbol{t}/\Pi_{n+1}(s)$, which yields

$$E[(X_{n+1} - X_n)^2|\mathcal{F}_n] = \frac{s^2}{\Pi_{n+1}^2(s)}\left[\frac{\boldsymbol{C}_n\boldsymbol{t}^2}{n+1} - \left(\frac{\boldsymbol{C}_n\boldsymbol{t}}{n+1}\right)^2\right], \tag{17}$$

and using the fact that each coordinate of $\boldsymbol{C}_n/(n+1)$ is bounded by 1 and Euler's formula for gamma function, the conditional second moment above is bounded by a constant multiple of $n^{-2s}$. Taking expectation and adding, we get $E[X_{n+1}^2]$ is bounded by a constant multiple of $\sum_0^n i^{-2s}$. This implies, $\{X_n\}$ is $L^2$-bounded if $s > 1/2$ and, $\{X_n/\sqrt{\log n}\}$ is $L^2$-bounded if $s = 1/2$. Thus, for $s > 1/2$, $X_n/\log n \to 0$ almost surely, as well as in $L^2$. For $s = 1/2$, $X_n/\log n \to 0$ in $L^2$.

We now show the convergence is almost sure also, when $s = 1/2$. For this, consider the random variables $Z_n = X_n/\log n$ and get

$$Z_{n+1} - Z_n = \frac{X_{n+1} - X_n}{\log(n+1)} - \frac{X_n}{\log n}\left[1 - \frac{\log n}{\log(n+1)}\right]. \tag{18}$$

Since $[1 - \frac{\log n}{\log(n+1)}]/\sqrt{\log n} \sim 1/n \log^{3/2} n$ and $X_n/\sqrt{\log n}$ is $L^2$-bounded (and hence $L^1$-bounded), we have $E[\sum_{k=2}^{n} \frac{|X_k|}{\sqrt{\log k}} \frac{1}{\sqrt{\log k}} \{1 - \frac{\log k}{\log(k+1)}\}]$ is bounded uniformly over $n$ and hence

$$\sum_{k=2}^{n} \frac{X_k}{\sqrt{\log k}} \frac{1}{\sqrt{\log k}} \left[ 1 - \frac{\log k}{\log(k+1)} \right]$$

converges absolutely almost surely. On the other hand, the first term of (18) is a martingale difference and, using (17) for $s = 1/2$, the conditional variance $E[(X_{n+1} - X_n)^2/\log^2(n+1)|\mathcal{F}_n]$ is bounded by a constant multiple of $[(n+1)\log^2(n+1)]^{-1}$, which is summable. Hence, the martingale $\{\sum_{k=1}^{n}(X_{k+1}-X_k)/\log(k+1)\}$ is $L^2$-bounded and thus converges almost surely. Combining the two observations above we get that $Z_n = X_n/\log n$ converges almost surely.

Thus $X_n/\log n$ converges to 0 almost surely, and in $L^2$, for all $s \geq 1/2$. Hence, from (16), we have

$$\frac{X_n}{\log n} = \frac{\boldsymbol{C}_n \boldsymbol{t}_2}{\log n \Pi_n(s)} - \frac{1}{\log n} \sum_{j=0}^{n-1} \frac{1}{j+1} \frac{\boldsymbol{C}_j \boldsymbol{t}_1}{\Pi_{j+1}(s)} \tag{19}$$

converges to 0 almost surely, as well as in $L^2$. But using (7) and Proposition 4.1(iii), we know that $\boldsymbol{C}_n \boldsymbol{t}_1/\Pi_n(s) \sim \Gamma(s+1)W_n/n^s \to \Gamma(s+1)V$ almost surely, as well as in $L^2$. Hence the second term in (19) converges to $\Gamma(s+1)V$ almost surely, as well as in $L^2$. Thus,

$$\frac{\boldsymbol{C}_n \boldsymbol{t}_2}{n^s \log n} \sim \frac{1}{\Gamma(s+1)} \frac{\boldsymbol{C}_n \boldsymbol{t}_2}{\Pi_n(s) \log n}$$

converges to $V$ almost surely, as well as in $L^2$. $\qquad\square$

**Remark 4.1.** As in the case of one dominant color, we have the correct scaling for all the linear combinations except when $\lambda = 0 < s$. (This situation arises only in the case of diagonalizable replacement matrix.) But $\boldsymbol{v}_2$ being an eigenvector of $R$ corresponding to $\lambda = 0$, we have $R\boldsymbol{v}_2 = \boldsymbol{0}$. Thus $\boldsymbol{C}_n \boldsymbol{v}_2 = \boldsymbol{C}_0 \boldsymbol{v}_2$ for all $n$.

The three-color urn model with two dominant colors can be easily extended to certain four-color models. We consider the reducible replacement matrix given in (5),

$$R = \begin{pmatrix} sQ & E \\ 0 & P \end{pmatrix},$$

where $P$ and $Q$ are $2 \times 2$ irreducible stochastic matrices, $0 < s < 1$. The eigenvalues of $Q$ are $\lambda$ and 1, with $|\lambda| < 1$. The eigenvalues of $P$ are $\beta$ and 1, with $|\beta| < 1$. Then $s\lambda$, $s$, $\beta$ and 1 are all eigenvalues of $R$. If $\boldsymbol{\xi}$ is an eigenvector of $Q$ corresponding to $\lambda$, then $\boldsymbol{v}_1 = (\boldsymbol{1}', \boldsymbol{0}')'$, $\boldsymbol{v}_2 = (\boldsymbol{\xi}', \boldsymbol{0}')'$ and $\boldsymbol{v}_4 = \boldsymbol{1}$ are eigenvectors of $R$ corresponding to $s$, $s\lambda$ and 1 respectively.

If $R$ is diagonalizable, then there is another eigenvector $\boldsymbol{v}_3$ corresponding to $\beta$. If $R$ is not diagonalizable, then one of its eigenvalues must repeat, namely $\beta$ must equal $s$ or $s\lambda$ and we denote the other by $\alpha$. In this case, we consider the Jordan decomposition $RT = TJ$, where $T$ is nonsingular. The fourth column $\boldsymbol{t}_4$ of $T$ can be chosen as $\boldsymbol{v}_4$. The first two columns $\boldsymbol{t}_1$ and $\boldsymbol{t}_2$ of $T$ can be chosen as the eigenvectors of $R$ corresponding to $\alpha$ and $\beta$. However, the third column $\boldsymbol{t}_3$ of $T$ will not be an eigenvector of $R$, yet the two-dimensional vector $\boldsymbol{\nu}$ formed by the lower half of $\boldsymbol{t}_3$ will be an eigenvector of $P$ corresponding to $\beta$. We shall only study $\boldsymbol{C}_n\boldsymbol{t}_3$ separately in the non-diagonalizable case.

The following three Propositions are suitable extensions of the three-color results of this section. The proofs are suitable modifications as well.

**Proposition 4.2.** *Consider a four-color urn model with the replacement matrix given by (5). Then:*

(i) $\boldsymbol{C}_n\boldsymbol{1}/(n+1) = 1$.

(ii) $\boldsymbol{C}_n/n \overset{\text{a.s.}}{\to} (0, 0, \boldsymbol{\pi}_P)$.

(iii) $(W_n, B_n)/n^s \overset{\text{a.s.}}{\to} \boldsymbol{\pi}_Q U$.

(iv) $\boldsymbol{C}_n\boldsymbol{v}_1/n^s \to U$ *almost surely, as well as in $L^2$.*

(v) *If $\lambda < 1/2$, then $\boldsymbol{C}_n\boldsymbol{v}_2/n^{s/2} \Rightarrow N(0, \; \frac{s^2\lambda^2}{s(1-2\lambda)}U\boldsymbol{\pi}_Q\boldsymbol{\xi}^2)$.*

(vi) *If $\lambda = 1/2$, then $\boldsymbol{C}_n\boldsymbol{v}_2/\sqrt{n^s \log n} \Rightarrow N(0, \; s^2\lambda^2 U\boldsymbol{\pi}_Q\boldsymbol{\xi}^2)$.*

(vii) *If $\lambda > 1/2$, then $\boldsymbol{C}_n\boldsymbol{v}_2/n^{s\lambda} \to V$ almost surely, as well as in $L^2$.*

*If we start with the initial vector $(W_0, B_0, G_0, Y_0)$, then $U$ and $V$ have the same distribution as the limit random variable in Theorem 3.1(iii) and the positive random variable in Theorem 3.1(vii), respectively, starting with initial vector $(W_0, B_0, G_0 + Y_0)$.*

Next we consider the linear combination $\boldsymbol{C}_n\boldsymbol{v}_3$ in the diagonalizable case.

**Proposition 4.3.** *In the four-color urn model with replacement matrix $R$ given by (5), assume that all the eigenvalues of $R$ are distinct. Then the following weak/strong laws hold for $\boldsymbol{C}_n\boldsymbol{v}_3$:*

(i) *If $\beta < 1/2$, then $\boldsymbol{C}_n\boldsymbol{v}_3/\sqrt{n} \Rightarrow N(0, \frac{\beta^2}{1-2\beta}\boldsymbol{\pi}_P\boldsymbol{\nu}^2)$.*

(ii) *If $\beta = 1/2$, then $\boldsymbol{C}_n\boldsymbol{v}_3/\sqrt{n \log n} \Rightarrow N(0, \; \beta^2\boldsymbol{\pi}_P\boldsymbol{\nu}^2)$.*

(iii) *If $\beta > 1/2$, then $\boldsymbol{C}_n\boldsymbol{v}_3/\Pi_n(\beta)$ is an $L^2$-bounded martingale and $\boldsymbol{C}_n\boldsymbol{v}_3/n^\beta$ converges almost surely to a non-degenerate random variable.*

Finally, we consider the case when the replacement matrix $R$ is not diagonalizable. As in the three-color urn model with a non-diagonalizable replacement matrix, the evolution of the linear combination $\boldsymbol{C}_n\boldsymbol{t}_3$ depends on the eigenvector of $R$ corresponding to the eigenvalue $\beta$. When $\beta < 1/2$, the effect of the contribution of the linear combination of this eigenvector is negligible. However, for $\beta \geq 1/2$, this provides the main contribution and the almost sure limit random variable depends on whether $\beta$ equals $s$ or $s\lambda$. To denote the limit random variable in a unified way, we define the random variable

$$W = \begin{cases} U, & \text{when } \beta = s, \\ V, & \text{when } \beta = s\lambda, \end{cases} \tag{20}$$

where $U$ and $V$ are the random variables defined in Proposition 4.2. Suppose $\beta \geq 1/2$. If $\beta = s$, then $\boldsymbol{t}_2 = \boldsymbol{v}_1$ and, by Proposition 4.2(iv), $\boldsymbol{C}_n \boldsymbol{t}_2/(n+1)^\beta = \boldsymbol{C}_n \boldsymbol{v}_1/(n+1)^s \to U = W$ almost surely, as well as in $L^2$. If $\beta = s\lambda$, then $\boldsymbol{t}_2 = \boldsymbol{v}_2$. Also $s < 1$ implies $\lambda > 1/2$. Hence by Proposition 4.2(vii), $\boldsymbol{C}_n \boldsymbol{t}_2/(n+1)^\beta = \boldsymbol{C}_n \boldsymbol{v}_2/(n+1)^{s\lambda} \to V = W$ almost surely, as well as in $L^2$. So, for non-diagonalizable $R$ and $\beta \geq 1/2$, we conclude $\boldsymbol{C}_n \boldsymbol{t}_2/(n+1)^\beta \to W$ almost surely, as well as in $L^2$.

**Proposition 4.4.** *Consider the four-color urn model with replacement matrix $R$ given by (5), where $R$ is not diagonalizable. Then, we have:*

(i) *If $\beta < 1/2$, then $\boldsymbol{C}_n \boldsymbol{t}_3/\sqrt{n} \Rightarrow N(0, \frac{\beta^2}{1-2\beta} \boldsymbol{\pi}_P \boldsymbol{\nu}^2)$.*

(ii) *If $\beta \geq 1/2$, then $\boldsymbol{C}_n \boldsymbol{t}_3/n^\beta \log n$ converges to $W$ almost surely, as well as in $L^2$, where $W$ is as defined in (20).*

**Remark 4.2.** As in two- and three-color urn models, we have correct scalings for all linear combinations except when $\lambda$ or $\beta$ becomes zero. If $\lambda = 0$, Proposition 4.2(v) gives $\boldsymbol{C}_n \boldsymbol{v}_2/n^{s/2} \xrightarrow{\text{P}} 0$. However, considering the three-color urn model $(W_n, B_n, G_n + Y_n)$, we get, from Remark 3.2, $\boldsymbol{C}_n \boldsymbol{v}_2 = \boldsymbol{C}_0 \boldsymbol{v}_2$.

In the case of the diagonalizable replacement matrix, if $\beta = 0$, we have a similar situation for the linear combination $\boldsymbol{C}_n \boldsymbol{v}_3$ in Proposition 4.3(i). However, as in Remark 4.1, $\boldsymbol{v}_3$ being an eigenvector of $R$ corresponding to $\beta = 0$, we have $R\boldsymbol{v}_3 = \boldsymbol{0}$ and $\boldsymbol{C}_n \boldsymbol{v}_3 = \boldsymbol{C}_0 \boldsymbol{v}_3$.

The situation becomes more interesting when $\beta = 0$ and the replacement matrix is not diagonalizable. Thus, we necessarily have $\beta = s$ or $\beta = \lambda s$, but $\beta$ being zero and $s$ being positive, only the second alternative is possible and further $\beta = \lambda = 0$. In this case, Proposition 4.4(i) gives $\boldsymbol{C}_n \boldsymbol{t}_3/\sqrt{n} \xrightarrow{\text{P}} 0$. The correct rate is given in the following proposition.

**Proposition 4.5.** *Consider the four-color urn model with the replacement matrix $R$ given by (5), which is not diagonalizable. Further assume the repeated eigenvalue of $R$ to be zero. Then*

$$\boldsymbol{C}_n \boldsymbol{t}_3/n^{s/2} \Rightarrow N(0, \boldsymbol{\pi}_Q \boldsymbol{\xi}^2 U/s),$$

*where $U$ is the limit random variable corresponding to $(W_n + B_n)/n^s$ obtained in Proposition 4.2(iv).*

**Proof.** Let $\boldsymbol{\chi}_n$ be the row vector as in the proof of Theorem 3.1(vii). Using $RT = TJ$ and the fact $\beta = 0$, we get $R\boldsymbol{t}_3 = \boldsymbol{t}_2 + \beta \boldsymbol{t}_3 = \boldsymbol{t}_2$. Hence using $\boldsymbol{t}_2 = (\boldsymbol{\xi}', \boldsymbol{0}')'$, the evolution equation for $C_n t_3$ is given by

$$C_{n+1} t_3 = C_n t_3 + \chi_{n+1} \xi.$$

The rest of the proof is similar to that of Theorem 4.2(i) and is omitted. □

# Acknowledgement

We thank the referee for an extremely careful reading of the manuscript and highly constructive comments.